\documentclass[10pt, reqno]{amsart}
\usepackage{amssymb,latexsym}
\usepackage{eucal}

\newtheorem*{proposition*}{}

\newtheorem{lemma}{Lemma}

\newtheorem*{Lem17.3}{Lemma 17.3}
\newtheorem*{Lem18.3}{Lemma 18.3}
\newtheorem*{Lem18.5}{Lemma 18.5}
\newtheorem*{Lem19.0}{Lemma 19.0}
\newtheorem*{Lem19.1}{Lemma 19.1}
\newtheorem*{Lem19.1.3}{Lemma 19.1.3}

\newtheorem*{Th}{Theorem}

\newcommand{\la}{{\langle}}
\newcommand{\ra}{{\rangle}}
\newcommand{\e}{{\varepsilon}}

\newcommand{\R}{\mathcal{R}}
\newcommand{\Ss}{\mathcal{S}}

\newcommand{\G}{\mathcal{F}}
\newcommand{\D}{\Delta}
\newcommand{\A}{\mathcal{A}}
\newcommand{\GG}{\mathcal{G}}

\newcommand{\X}{\mathcal{X}}

\newcommand{\p}{\partial}

\newcommand{\C}{\mathcal{C}}

\newcommand{\wtl}{\widetilde}

\begin{document}

\title[On identities in groups of fractions of  cancellative semigroups]
{On identities in groups of fractions of  cancellative
semigroups}
\author{S.V. Ivanov}
\address{Department of Mathematics\\
University of Illinois \\
Urbana,  IL 61801, USA}
\email{ivanov@math.uiuc.edu}
\thanks{The first author is supported in part by NSF grant DMS 00-99612}

\author{A.M. Storozhev}
\address{Australian Mathematics Trust\\
University of Canberra \\
Belconnen, ACT 2616, Australia}
\email{andreis@amt.canberra.edu.au }


\begin{abstract}
To solve two problems of Bergman stated in 1981, we construct a
group $G$ such that $G$ contains a free noncyclic subgroup (hence,
$G$ satisfies no group identity) and $G$, as a group, is generated
by its subsemigroup that satisfies a nontrivial semigroup
identity.
\end{abstract}
\maketitle

\section{Introduction}

A semigroup $S$ is called {\em cancellative} if for arbitrary
$a, b, x \in S$ either of equalities $xa=xb$, $ax=bx$ implies
that  $a=b$.

A semigroup $S$ is said to satisfy {\em left} (resp. {\em
right}) {\em Ore's condition} if  for arbitrary $a, b \in S$
there are $x, y \in S$ such that $x a = y b$ (resp. $ax = by$).

It is well known that a cancellative semigroup $S$ with left (or
right) Ore's condition embeds in a group $\G(S)= S^{-1} S = S
S^{-1}$ of its fractions (see \cite{M53}, \cite{NT63},
\cite[Theorem 1.23]{CP67}; recall that, in general, a cancellative
semigroup need not embed in a group, see \cite{M37}) and a group
$\G(S)$ of fractions of $S$ is unique in the sense that if $G$ is
a group such that $G$ contains $S$ and $G$, as a group,  is
generated by $S$ then $G$ is naturally isomorphic to $\G(S)$.

In particular, if $S$ is a cancellative semigroup with a
nontrivial identity then $S$ satisfies (both) Ore's conditions
and so $S$ embeds in the  group $\G(S)$ of its fractions.

If $S$ is an abelian cancellative semigroup then the group of
its fractions is also abelian. Mal'tsev \cite{M53} proved that
if $S$ is a nilpotent of class $n$ semigroup (that is, $S$ is
given by a semigroup identity of the form $X_n \equiv Y_n$,
where $X_n, Y_n$ are words in letters $x, y, u_1, \dots, u_n$,
inductively defined by $X_0 = x$, $Y_0 = y$ and, for $n \ge 0$,
by $X_{n+1} = X_n u_{n+1} Y_n$, $Y_{n+1} = Y_n u_{n+1} X_n $),
then $\G(S)$ is a nilpotent group of class $n$. It was also
shown in \cite{M53}  that a group $K$ is  nilpotent of class
$n$ if and only if the identity $X_n \equiv Y_n$ holds in $K$.

In connection with the results stated above, Bergman \cite{B81}
posed the following natural problem (which, in a different form,
is also mentioned by Shevrin and Suhanov in \cite[Question
11.1]{SS89}): Suppose a cancellative semigroup $S$ satisfies a
nontrivial semigroup identity. Does then the group $\G(S)$ of
fractions of $S$ satisfy all the semigroup identities that hold in
$S$? Bergman \cite{B81} also proposed a weaker problem: If a
cancellative semigroup $S$ satisfies a nontrivial semigroup
identity,  must then the group $\G(S)$ of its fractions satisfy a
nontrivial group identity?

Krempa and Macedo\'nska \cite{KM92} proved that if an
 identity of the form $x^ny^n \equiv y^nx^n$, $n \ge 1$, holds  in a
cancellative semigroup $S$ then every semigroup identity of $S$
holds in  $\G(S)$ (hence, both Bergman's problems have positive
solutions for the identity $x^ny^n \equiv y^nx^n$). Krempa and
Macedo\'nska \cite{KM92} also observed that the same result holds
for the Mal'tsev nilpotent identity $X_n \equiv Y_n$, $n \ge 1$.

In this paper, we construct an $m$-generator group $\GG = \la \bar
a_1, \dots, \bar a_m \ra$, $m >1$, generated by $\bar a_1, \dots,
\bar a_m$, such that $\GG$ contains a free noncyclic subgroup
(hence, no group identity holds in $\GG$) and the subsemigroup
$\Ss =  \la \bar a_1, \dots, \bar a_m \ra^+$ of $\GG$, generated
by $\bar a_1, \dots, \bar a_m$, satisfies a nontrivial semigroup
identity (namely, $w_L(x,y) \equiv w_R(x,y)$, see
\eqref{wl}--\eqref{wr} below). Since $\Ss$ is obviously
cancellative and $\GG$, as a group,  is generated by $\Ss$, it
follows that both Bergman's problems are solved in the negative.

It should be pointed out that, in the 1980's,  Rips suggested a
similar approach to solve Bergman's problems in the negative.
However, to the best of our knowledge, no details of his
construction are available.

Now we will give details of our construction. Consider the
following semigroup words
\begin{equation} \label{wl}
w_L(x,y) = (x^dy^d)^{n+1^2} x (x^dy^d)^{n+2^2} x \dots
(x^dy^d)^{n+(h_0-1)^2} x (x^dy^d)^{n+h_0^2} x
\end{equation}
and
\begin{multline} \label{wr}
w_R(x,y) = x (x^dy^d)^{n-h_0^2(2h_0-3)^2} x
(x^dy^d)^{n+(2h_0-1)^2} x (x^dy^d)^{n+(2h_0-2)^2} \dots \\
\dots x (x^dy^d)^{n + (h_0+2)^2} x (x^dy^d)^{n+(h_0+1)^2} ,
\end{multline}
where  $h_0 = h/2, d, n$ are sufficiently large positive integers
(with $n \gg d \gg h \gg 1$). Observe that it follows from the
equality $1^2 +2^2 + \dots + k^2 = \frac {k(2k+1)(k+1)}{6}$ that
the sum of exponents on $x$ (resp. $y$) in $w_L(x,y)$ is equal to
that on $x$ (resp. $y$) in $w_R(x,y)$. Note that the word
\begin{equation} \label{w}
w(x, y) = w_L(x,y) w_R(x,y)^{-1} ,
\end{equation}
was introduced in \cite{OS96} to construct the first example of a
finitely generated group $G_{0}$ such that $G_{0}$ satisfies a
nontrivial semigroup identity (namely, $w_L(x,y) \equiv w_R(x,y)$)
and $G_{0}$ is not a periodic extension of a locally nilpotent
group. However, the parameter $h$, used in the construction of
$w(x,y)$ in \cite{OS96}, should actually read $h/2$ and so we put
$h_0 = h/2$ in \eqref{wl}, \eqref{wr} (recall that, according to
the condition $R$ of \cite[Sects. 25]{O89}, every 2-cell of second
type must have exactly $h$ long sections).

We also recall that the question on the existence of a finitely
generated group that satisfies a nontrivial semigroup identity and
is not a periodic extension of a locally nilpotent group was asked
by Longobardi, Maj and Rhemtulla \cite{LMR95} and goes back to a
result of Rosenblatt \cite{R74} (rediscovered in \cite{LMR95})
that a finitely generated solvable group which contains no free
subsemigroups of rank 2 is nilpotent-by-finite. Here is our main
result.

\begin{Th}
There exists an $m$-generator torsion-free group $\GG= \la \bar
a_1,$ $\dots, \bar a_m \ra$, $m \ge 2$,  such that $\GG$ contains
a subgroup isomorphic to a free group of rank $2$ $($e.g.,
generated by $\bar a_1 \bar a_2^{-1}\bar a_1$ and $\bar a_2 \bar
a_1^{-1}\bar a_2$$)$ and the nontrivial semigroup identity
$w_L(x,y) \equiv w_R(x,y)$ $($see \eqref{wl}--\eqref{wr}$)$ holds
in the subsemigroup $\Ss = \la \bar a_1, \dots, \bar a_m \ra^+$ of
$\GG$, generated by $\bar a_1, \dots, \bar a_m$.
\end{Th}

\section{Inductive Construction}

To prove this Theorem, we will inductively construct a
presentation for the group  $\GG$ by means of generators and
defining relations and use the geometric machinery of graded
diagrams,  developed by Ol'shanskii \cite{O85}, \cite{O89} and
refined by the authors in \cite{IS03}, to study this group. In
particular, we will use the notation and terminology of \cite{O89}
and all notions that are not defined in this paper can be found in
\cite{O89}.

As in \cite{O89},  we will use numerical parameters
$$
\alpha \succ \beta  \succ \gamma  \succ \delta \succ \e \succ
\zeta \succ \eta \succ \iota
$$
and $h = \delta^{-1}$, $d = \eta^{-1}$, $n = \iota^{-1}$ ($h,
d, n$ were already used in \eqref{wl}--\eqref{wr}) and employ
the least parameter principle (LPP) (according to LPP a small
positive value for, say,  $\zeta$ is chosen to satisfy all
inequalities whose smallest (in terms of the relation $\succ$)
parameter is $\zeta$).

Let $\A = \{ a_1, \dots, a_m \}$  be an alphabet, $m > 1$, and
$F(\A)$ be the free group in $\A$. Elements of $F(\A)$ are
referred to as words in $\A^{\pm 1} = \A \cup \A^{-1}$ or just
words. We will say that $U \in F(\A)$ is a {\em positive} word
if $U$ has no occurrences of $a_1^{-1}, \dots, a_m^{-1}$. A
word $V$ will be called {\em regular} if at least one of $V$,
$V^{-1}$ contains a positive subword of length 3.

Denote $G(0) = F(\A)$ and let the set $\R_0$ of defining words
of rank 0  be empty. To define the group $G(i)$ by induction on
$i \ge 1$, assume that the group $G(i-1)$ is already
constructed by its presentation
$$
G(i-1) = \la \A \; \| \; R=1, R \in \R_{i-1} \ra .
$$
Let $\X_i$ be a  set of  words (in $\A^{\pm 1}$) of length $i$,
called {\em periods of rank $i$}, which is maximal with respect
to the following two properties:

\begin{enumerate}
\item[(A1)] If $A \in \X_i$ then $A$ (that is, the image of $A$
in $G(i-1)$) is not conjugate in $G(i-1)$ to a power of a word
of length $< |A| = i$.

\item[(A2)] If $A$, $B$ are distinct elements of $\X_i$ then
$A$ is not conjugate in $G(i-1)$ to $B$ or $B^{-1}$.
\end{enumerate}

If the images of two words $X$, $Y$ are equal in the group
$G(i-1)$, $i \ge 1$, then we will say that $X$  is {\em equal
in rank} $i-1$ to $Y$ and write $ X \overset {i-1} = Y$.
Analogously, we will say that two words $X$, $Y$ are {\em
conjugate in rank} $i-1$ if their images are conjugate in the
group $G(i-1)$. As in  \cite{O89}, a word $A$ is called {\em
simple} in rank $i-1$, $i \ge 1$, if $A$ is conjugate in rank
$i-1$ neither to a power $B^\ell$, where $|B| <|A|$ nor to a
power of period of some rank $\le i-1$. We will also say that
two pairs $(X_1, X_2)$, $(Y_1, Y_2)$ of words are conjugate in
rank $i-1$, $i \ge 1$, if there is a word $W$ such that  $X_1
\overset {i-1} = W Y_1 W^{-1}$ and $X_2 \overset {i-1} = W Y_2
W^{-1}$.

Consider the set of all possible pairs $(X, Y)$ of words in
$\A^{\pm 1}$. This set is partitioned by equivalence classes
$\C_\ell$, $\ell=1,2,\dots$,  of the equivalence relation
$\sim$ defined by $(X_1, Y_1) \sim (X_2, Y_2)$ if and only if
the pairs  $( X_1^d Y_1^d, w(X_1, Y_1) )$ and $( X_2^d Y_2^d,
w(X_2, Y_2))$, where $w(x,y)$ is given by \eqref{w}, are
conjugate in rank $i-1$. It is convenient to enumerate (in some
way)
$$
\C_{A^f, 1}, \C_{A^f, 2}, \dots
$$
 all classes of pairs $(X, Y)$ such that $w(X, Y)
\overset {i-1} \neq 1$ and  $X^d Y^d$ is conjugate in rank
$i-1$ to some power $A^f$, where $A \in \X_i$ and $f$  are
fixed.

It follows from definitions that every class $\C_{A^f, j}$  contains
a pair
$$
(X_{A^f, j}, \bar Y_{A^f, j} )
$$
with the following properties. The word $X_{A^f, j}$ is
graphically (that is, letter-by-letter) equal to a power of
$B_{A^f, j}$, where $ B_{A^f, j}$ is simple in rank $i-1$ or a
period of rank $\le i-1$; $\bar Y_{A^f, j} \equiv Z_{A^f, j}
Y_{A^f, j} Z_{A^f, j}^{-1}$, where the symbol \ '$\equiv$' \
means the graphical equality, $Y_{A^f, j}$ is graphically equal
to a power of $C_{A^f, j}$, where $C_{A^f, j}$ is simple in
rank $i-1$ or a period of rank $\le i-1$. We can also assume
that if $D_1 \in \{ A, B_{A^f, j}, C_{A^f, j} \}$  is conjugate
in rank $i-1$ to $D_2^{\pm 1}$, where $D_2 \in \{ A, B_{A^f,
j}, C_{A^f, j} \}$, then $D_1 \equiv D_2$. Finally, the word
$Z_{A^f, j}$ is picked for fixed $X_{A^f, j}$, $Y_{A^f, j}$ so
that the length $|Z_{A^f, j}|$ is minimal (and the pair
$(X_{A^f, j}, Z_{A^f, j} Y_{A^f, j} Z_{A^f, j}^{-1})$ belongs
to $\C_{A^f, j}$). Similar to \cite{O85}, \cite{O89},
\cite{IS03}, the triple $(X_{A^f, j}, Y_{A^f, j}, Z_{A^f, j} )$
is called an $(A^f, j)$-{\em triple} corresponding to the class
$\C_{A^f, j}$ (in rank $i-1$).

Now for every class $\C_{A^f, j}$, where $A$ is a regular word,
we pick a corresponding  $(A^f, j)$-{triple}
$$
(X_{A^f, j},  Y_{A^f, j}, Z_{A^f, j} )
$$
in rank $i-1$ and construct a defining word $R_{A^f, j}$ of rank $i$ as
follows.
Pick a word $W_{A^f, j}$ of minimal length so that
$$
X_{A^f, j}^d  \bar Y_{A^f, j}^d \overset {i-1} = W_{A^f, j}
A^{f} W_{A^f, j}^{-1} .
$$
Let   $T_{A^f, j}$ be a word of minimal length such that
\begin{gather*}
T_{A^f, j}  \overset {i-1}  = W_{A^f, j}^{-1}  X_{A^f, j}
W_{A^f, j}.
\end{gather*}
According to  \eqref{w} (see also \eqref{wl}--\eqref{wr}), we
set
\begin{multline}
R_{A^f, j}  =  A^{(n+1^2)f} T_{A^f, j} A^{(n+2^2)f} T_{A^f, j}
\dots A^{(n+(h_0-1)^2)f} T_{A^f, j} A^{(n+h_0^2)f} T_{A^f, j} \\
 A^{(-n-(h_0 + 1)^2)f} T_{A^f, j}^{-1}
A^{(-n-(h_0 + 2)^2)f} T_{A^f, j}^{-1}
 \dots
A^{(-n-(2h_0-2)^2)f} T_{A^f, j}^{-1} \\  A^{(-n-(2h_0-1)^2)f}
T_{A^f, j}^{-1}A^{(-n+h_0^2(2h_0-3)^2)f} T_{A^f, j}^{-1}.
\label{D:3}
\end{multline}
It follows from definitions that the word $R_{A^f, j}$ is
conjugate in rank $i-1$ (by $W_{A^f, j}^{-1}$) to the word
$w(X_{A^f, j}, \bar Y_{A^f, j}) \overset {i-1} \neq 1$.

The set $\Ss_i$ of defining words of rank $i$ consists of all
possible words $R_{A^f, j}$  given by (\ref{D:3}) (over all
equivalence classes $\C_{A^f, j}$, where $A \in \X_i$ is a
regular word). Finally, we put $\R_i = \R_{i-1} \cup \Ss_i$ and
set
\begin{gather}
G(i) = \la \A \; \| \; R=1, R \in  \R_i \ra .  \label{D:5}
\end{gather}
The inductive definition of  groups $G(i)$, $i \ge 0$, is now
complete and we can consider the limit group $\GG$ given by
defining words of all ranks $j =1,2, \dots$
\begin{gather}
 \GG
= G(\infty) = \la \A \; \| \; R=1, R \in \cup_{j=0}^\infty \R_j
\ra . \label{D:6}
\end{gather}

\section{Proof of Theorem}

First we will establish several lemmas.

\begin{lemma} \label{L1}
The presentation $(\ref{D:5})$ of $G(i)$ satisfies the  condition
$R$ of \cite[Sect. 25]{O89}.
\end{lemma}

\begin{proof}
Recall that in article  \cite{OS96} the free group
$F(\A)/w(F(\A))$ of the variety of groups, defined by the identity
$w(x,y) \equiv 1$, see \eqref{w}, is constructed. The inductive
construction of a presentation for the group $F(\A)/w(F(\A))$ in
\cite{OS96} is analogous to our construction of the group $ \GG =
G(\infty)$ (the only difference is that $A \in \X_i$ is arbitrary
in \cite{OS96} and now we require that $A$ be regular) and many
technical results of \cite{OS96} are reproved in our situation
without any changes. In particular, repeating arguments of proofs
of Lemmas 4--6 \cite{OS96}, we can prove Lemma \ref{L1}.
\end{proof}

Now suppose that $X, Y$ are some words such that  $[X, Y] =
XYX^{-1} Y^{-1}  \overset {i} \neq 1.$ Conjugating the pair
$(X, Y)$ in rank $i$, if necessary, we can assume that $X
\equiv B^{f_B}$,  $Y \equiv Z C^{f_C} Z^{-1}$, where each of
$B, C$ is either simple in rank $i$ or a period of some rank
$\le i$ and, when $B^{f_B}$, $C^{f_C}$ are fixed, the word $Z$
is picked to have minimal length. We also consider the
following equality
\begin{gather*}
X^d Y^d \overset {i} = W_A A^{f_A} W_A^{-1}  ,
\end{gather*}
where $ A$ is either simple in rank $i$ or a period of some
rank $\le i$ and the conjugating word $W_A$ is picked (when $A$
is fixed) to have minimal length. Without loss of generality,
we can also assume that if $D_1 \in \{ A, B, C \}$ is conjugate
in rank $i$ to $D_2^{\pm 1}$, where $D_2 \in \{ A, B, C\}$,
then $D_1 \equiv D_2$.

\begin{lemma}  \label{L2}
In the foregoing notation, the following inequalities hold
\begin{gather*}
0 < | f_A | \le 100 \zeta^{-1} ,  \quad  \max (| B^{d f_B} |, |
C^{d f_C} | ) \le   \zeta^{-2} | A^{f_A} | ,  \quad  |Z| < 3
\zeta^{-2} | A^{f_A} | .
\end{gather*}
\end{lemma}

\begin{proof} These inequalities are proved analogously  to
inequalities (13)--(14) of Lemma~6 \cite{IS03} (see also Lemma
3 \cite{S94}).
\end{proof}

\begin{lemma}  \label{L3}
 Let $\wtl  X$ and $\wtl Y$ be positive words.
 Then $w(\wtl X,\wtl Y)=1$ in the group $G(\infty)$ given by presentation
$(\ref{D:6})$.
\end{lemma}

\begin{proof} Arguing on the contrary, assume that
\begin{equation}
w(\wtl X,\wtl Y) \ne 1 \label{L2:0}
\end{equation}
in $G(\infty)$. Let $A$ be a period of some rank such that
$A^{f_A}$ for some $f_A$ is conjugate in $G(\infty)$ to $\wtl
X^d \wtl Y^d$. (The existence of such an $A$ follows from
definitions; see also Lemma 18.1 \cite{O89}.) Note that, in
view of \eqref{L2:0}, $[\wtl X, \wtl Y ]   \neq 1$ in
$G(\infty)$. Hence, by Lemma \ref{L2}, we can replace the pair
$(\wtl X, \wtl Y )$ by a conjugate in the group $G(\infty)$
pair $(X, Y)$  such that $X \equiv B^{f_B}$, $Y \equiv Z
C^{f_C} Z^{-1}$, and $X^d Y^d = W_A A^{f_A} W_A^{-1}$ in
$G(\infty)$, where $B$, $C$ are some periods and inequalities
of Lemma \ref{L2} hold. (Note that $X, Y$ need not be positive
words.) In particular,
\begin{equation} \label{est}
|X^d Y^d | \le | B^{d f_B} | + |C^{d f_C} | +2 |Z| <
 8 \zeta^{-2}
|A^{f_A}| .
\end{equation}

Let $\wtl \D$ be a reduced annular diagram of rank $\wtl i$ for
conjugacy of $\wtl X^d \wtl Y^d$ and $A^{f_A}$. By Theorem 22.4
\cite{O89} applied to $\D$, we have
\begin{equation}
(1- \beta) |A^{f_A}| \le |\wtl X^d\wtl Y^d|. \label{L2:1}
\end{equation}

It follows from the definition of the word $w(x,y)$ (see
\eqref{w}) that the sums of exponents on $x$ and $y$ in
$w(x,y)$ are 0.  Since all defining words of presentation
\eqref{D:6} are values of $w(x,y)$, these defining  words are
in the commutator subgroup of $F(\A)$. Hence, the total sum
$\sigma(\wtl X^d \wtl Y^d) + \sigma(A^{-f_A})$ of exponents on
letters $a_1, \dots, a_m$ in words $\wtl X^d \wtl Y^d$ and
$A^{-f_A}$ is equal to 0. Let $\sigma^+(A^{f_A})$ denote the
sum of all positive exponents on letters $a_1, \dots, a_m$ in
$A^{f_A}$. Since $\wtl X^d \wtl Y^d$ is conjugate in
$G(\infty)$ to $A^{f_A}$ and $\wtl X^d, \wtl Y^d$ are positive
words, it follows that
\begin{equation}
\sigma^+(A^{f_A})  \ge |\wtl X^d \wtl Y^d|. \label{L2:1.1}
\end{equation}

In view of inequalities \eqref{L2:1}--\eqref{L2:1.1}, we
further have
\begin{equation}
 \sigma^+(A^{f_A})  \ge (1- \beta) |A^{f_A}| . \label{L2:1.2}
\end{equation}

Assume that $A$ is not regular. Then the sum of positive
exponents on $a_1, \ldots, a_m$ in $A^{f_A / |f_A|}$ is at most
$\frac 23 (|A|+1)$. Therefore,  $ \sigma^+(A^{f_A})  \le \frac
23 |f_A| (|A|+1)$ and, in view of \eqref{L2:1.2}, we get $\frac
23 (|A|+1) \ge  (1-\beta) |A|$ which implies that $|A|<3$
(LPP). On the other hand, it follows from Lemma \ref{L2} that
$$
|A|\ge |f_A|^{-1} \zeta^2 |X^d| \ge 100^{-1} \zeta^3 d > 3
$$
(LPP: $\zeta \succ \eta=d^{-1}$). This contradiction to $|A| <
3$ shows that $A$ is regular.

Now we consider a reduced annular diagram $\D$ of some rank
$i'$ for conjugacy of $X^d Y^d$ and $A^{f_A}$. By Lemmas
\ref{L1} and 22.1 \cite{O89}, $\D$ can be cut into a simply
connected diagram $\D_1$ along a simple path $t$ which connects
points on distinct components of $\p \D$ with $| t| < \gamma
|\p \D |$. Therefore, in view of estimate \eqref{est},
\begin{multline*}
|\p \D_1 | < (1 + 2 \gamma)|\p \D | = (1 + 2 \gamma) (|X^d Y^d|
+|A^{f_A}|) <\\ < (1 + 2 \gamma)( 8\zeta^{-2}+ 1) |A^{f_A}| <
10^3\zeta^{-3} |A| < \tfrac 12 n |A|
\end{multline*} (LPP: $\gamma \succ \zeta \succ \iota = n^{-1}$).
Then, by Lemmas \ref{L1}, 20.4 and 23.16 \cite{O89} applied to
$\D_1$, the diagram $\D_1$ contains no 2-cells of rank $ \ge
|A|$, whence  $\D_1$, $\D$ are diagrams of rank $|A| -1$. Since
$A \in \X_{|A|}$ is a regular word, it follows from the
construction of defining words of rank $|A|$ that there will be
a defining word in $\Ss_{|A|}$ which guarantees that $w(X, Y)
\overset {|A|} = 1$. Consequently, $w(\wtl X, \wtl Y) = 1$ in
$G(\infty)$ and a contradiction to assumption (\ref{L2:0})
proves Lemma \ref{L3}.
\end{proof}

{\em Proof of Theorem.} It follows from definitions and Lemmas
\ref{L1}, 25.2 \cite{O89} that the  group $\GG =G(\infty) = \la
a_1, \dots, a_m \ra$, given by presentation (\ref{D:6}),  is
torsion-free. By Lemma \ref{L3}, the subsemigroup $\Ss = \la a_1,
\dots, a_m \ra^+$ of $\GG$, generated by $a_1, \dots, a_m$,
satisfies the nontrivial semigroup  identity $w_L(x,y)\equiv
w_R(x,y)$. It remains to show that $\GG$ contains a subgroup
isomorphic to a free group of rank 2. To do this, consider words
$V_1 =a_1 a_2^{-1} a_1$ and $V_2=a_2 a_1^{-1} a_2$.  Assume that
there is a nonempty cyclically reduced word $U(V_1, V_2)$ in
$V_1^{\pm 1}$, $V_2^{\pm 1}$ such that $U(V_1, V_2) =1$ in the
group $\GG =G(\infty)$. Then we can consider a reduced disk
diagram $\D$ of positive rank the label of whose boundary $\p \D$
is $U(V_1, V_2)$.  It follows from the definition of defining
words of  $\GG =G(\infty)$ and Theorem 22.2 \cite{O89} applied to
$\D$ that $U(V_1, V_2)$ is a regular word. This, however, is
obviously false. This contradiction shows that the subgroup $\la
V_1, V_2 \ra$ of $\GG$ is free of rank 2 and Theorem is proved.
\qed

\end{document}